\documentclass [12pt] {article}
\title {A link polynomial via a vertex-edge-face state model}
\author{Thomas Fiedler}
\usepackage{amssymb,amsfonts,epsfig}
\begin{document}
\newtheorem{proposition}{Proposition}
\newtheorem{theorem}{Theorem}
\newtheorem{lemma}{Lemma}
\newtheorem{corollary}{Corollary}
\newtheorem{example}{Example}
\newtheorem{remark}{Remark}
\newtheorem{definition}{Definition}
\newtheorem{question}{Question}
\newtheorem{conjecture}{Conjecture}
\maketitle
\begin{abstract}
We construct a  2-variable link polynomial, called  $W_L$,  for classical links  by considering simultaneously the Kauffman state 
models for the Alexander and for the Jones polynomials. We conjecture that this polynomial is  the product of two 
1-variable  polynomials, one of which is the Alexander polynomial.

We refine $W_L$ to an ordered set of 3-variable polynomials for those links in 3-space which contain a Hopf link as a 
sublink. 

\footnote{
classical links, quantum invariants, state models}
\end{abstract}

\section{Introduction and results}
We work in the smooth category. All 2- and 3-dimensional manifolds are oriented.

The present paper should have been written 20 years ago, when the state models for the Alexander and for the Jones
polynomials were discovered. 

Let us fix a coordinate system $(x ,y ,z)$ in $\mathbb{R}^3$ and let $pr : (x ,y ,z) \to (x ,y)$ be the standard projection.
Let $S^3 = \mathbb{R}^3 \cup \infty$ and let $S^2 = (x ,y)-plane \cup \infty$. Let $L \hookrightarrow \mathbb{R}^3$ be an
 oriented link. We represent links as usual by diagrams D
with respect to $pr$ (see e.g. \cite{BZ}). Let $A = z-axes \cup \infty$. A is a {\em meridian} of D if there is a half-plane 
bounded by the z-axes 
which intersects D in exactly one point. If the intersection index at this point is $+1$ then A is called a {\em positive}
meridian. Let D and D' represent the same knot in $S^3$ and such that A is a positive meridian 
for both of them.
Then, as well known, D and D' represent the same knot in the solid torus $\mathbb{R}^3 \setminus z-axes$
(see e.g. \cite{F1}). For links 
we have to replace isotopy by isotopy which preserves the distinguished components (for which we have choosen the 
meridians).

More generally, in this paper we will study oriented links in the solid torus $V = \mathbb{R}^3 \setminus z-axes$.

Let $D$ be an oriented link diagram. Its projection $pr(D)$ is an oriented graph in the annulus. We call the double points 
in $pr(D)$ the {\em vertices}, the arcs which connect the double points the {\em edges} and the components of its 
complement in the plane or
in the annulus the {\em faces}.

Kauffman has constructed a vertex-face  state model for the Alexander polynomial (see \cite{K1}) and a vertex-edge state
model for the Jones polynomial (see \cite{K2}). For the convenience of the reader we remaind the definitions here.

Let D be an oriented and connected diagram of a link in $\mathbb{R}^3$. There are exactly two more regions in 
$\mathbb{R}^2 \setminus pr(D)$ than crossings of D. We mark two adjacent regions by stars (i.e. their boundaries in 
$\mathbb{R}^2$
have a common edge in the 4-valent graph $pr(D)$ ). A state T assignes now to each crossing of $pr(D)$ a dot in 
exactly one of the four local quadrants in $\mathbb{R}^2 \setminus pr(D)$ and such 
that in each region of $\mathbb{R}^2 \setminus pr(D)$, besides the regions marked by the stars , there is exactly one 
dot.
To each quadrant we associate a monomial as shown in Fig. 1. Notice, that if we switch the crossing then the 
monomial is replaced by its inverse.

\begin{figure}
\centering 
\psfig{file=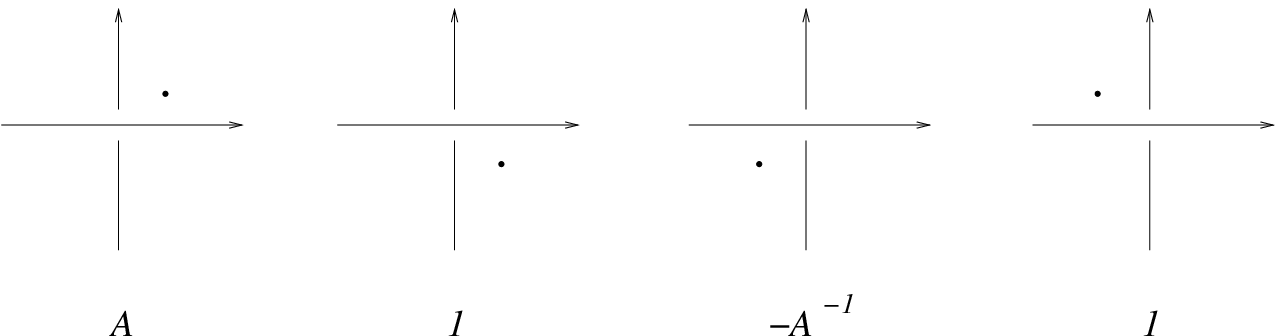}
\caption{}
\end{figure}

To each state T we associate now the product of the monomials corresponding to the dots. The {\em Alexander 
polynomial\/}
 $\Delta_L(t) \in \mathbb{Z}[t ,t^{-1}]$ is then the sum of all products of monomials over all possible states T. Here , 
L denotes the link represented by D.
Kauffman shows that $\Delta_L(t)$ is an isotopy invariant of L , that it does not depend on the choice of the adjacent 
regions with a star and that it actually coincides with the Alexander polynomial.

The {\em Kauffman bracket} in $\mathbb{R}^3$ is defined as follows: a state S  splitts $pr(D)$ at each double point in 
exactly one of the 
two possible ways. To each such splitting we assigne a monomial as shown in Fig. 2.

\begin{figure}
\centering 
\psfig{file=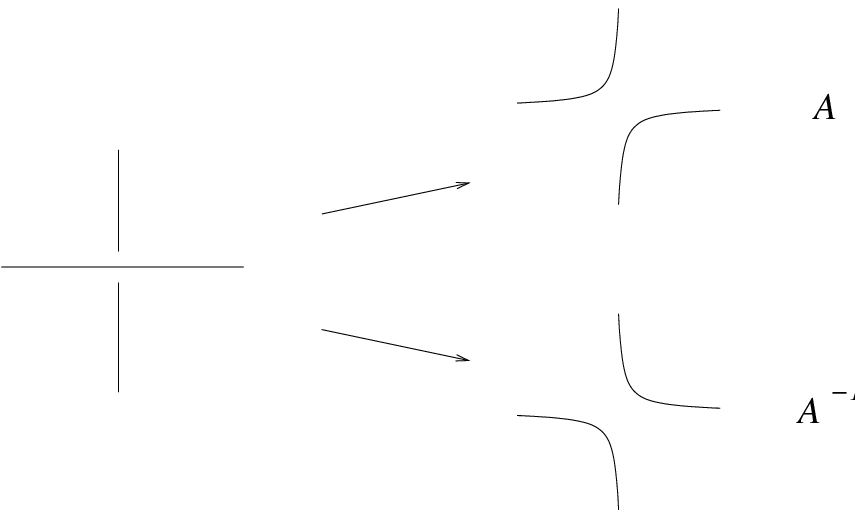}
\caption{}
\end{figure}

 Notice, that again if we switch the
crossing then the monomial is replaced by its inverse. Let $/S/$ be the number of circles which are the result of splitting
all double points of $pr(D)$ according to the state S. To each state S we associate the polynomial $<D,S>$  which is the 
product of all the monomials coming from the double points with $d^{/S/}$, where $d = -A^2-A^{-2}$. The Kauffman 
bracket $<D>$ is then the sum of the polynomials $<D,S>$ over all states S. The Kauffman bracket is invariant under
Reidemeister moves of type II and III but it is not invariant under Reidemeister moves of type I. But Kauffman shows 
that $(-A)^{-3w(D)}<D>$ , where $w(D)$ is the {\em writhe} (compare e.g.  \cite{BZ}), is a link invariant which coincides (up to normalisation for the unknot)
with the {\em Jones polynomial}  $V_L(t)$ for $A = t^{-1/4}$.

In the case of the solid torus this state model can be refined (see \cite{HP}): let $/S/$ be now only the number of contractible 
circles in the annulus and let $[S]$ be the number of non contractible circles. We replace then $d^{/S/}$ by 
$d^{/S/}h^{[S]}$, where $h$ is a new independent variable. For the Conway and the Kauffman skein modules of the solid
 torus see \cite{T}.

The starting point of the present paper is the following simple observation (already used in \cite{F3}):  let D be a connected 
diagram in the solid torus 
$V = \mathbb{R}^3 \setminus z-axes$ and
such that the corresponding link L is not contained in a 3-ball in the solid torus. (Evidently, we can always make D 
connected by performing just some Reidemeister II moves.) Then there are two {\em canonical 
regions\/}
in $S^2 \setminus pr(D)$, namely those which contains $\infty$ and those which contains the origin $(0 ,0)$. We mark 
the canonical regions by the stars.
Notice, that our star-regions are adjacent if and only if A is a meridian up to isotopy of the diagram D. 

The important point is that we have no longer to prove invariance under the choice of the stars. This gives us the 
possibility to consider both types of state sums simultaneously.

\begin{definition}
Let T be a Kauffman state for the Alexander polynomial. At each crossing we consider the two possible splittings 
(indicated by a small dash in the figures), i.e. the two Kauffman states S for the Jones polynomial. 
We associate to each positive crossing for each T and each S a (complex) monomial as shown in Fig. 3. Here, $x$  
$y$ and $z$ are
independent variables.

\begin{figure}
\centering 
\psfig{file=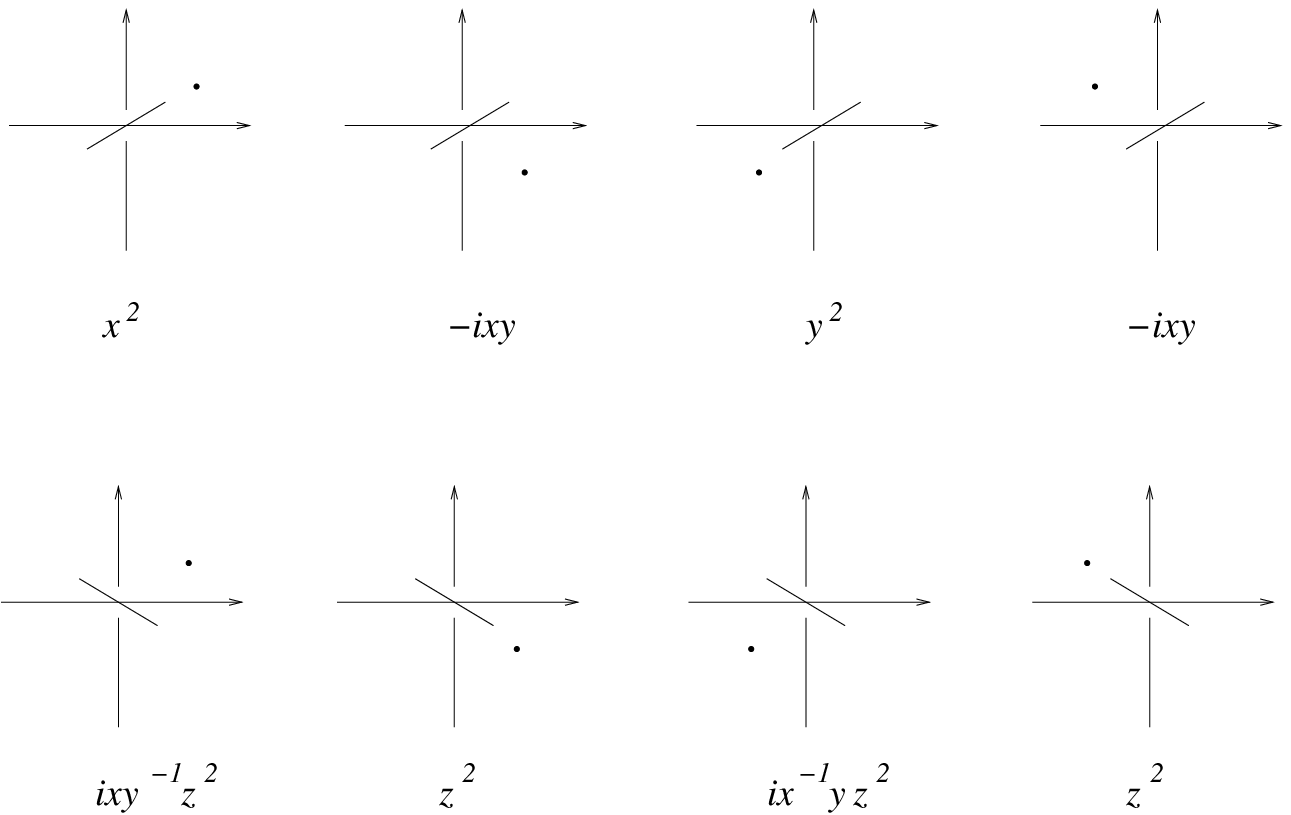}
\caption{}
\end{figure}

 If we switch the crossing to a negative one then we associate for the same T and S the
inverse monomial.

We set $d = ixyz^{-2} - ix^{-1}y^{-1}z^2$ and $h$ is an independent variable as previously.

To each couple of states T and S , called a {\em double state} $(T,S)$, we associate the {\em double bracket} $<D,T,S>$ which is the product of the above
defined monomials over all crossings of D and $d^{/S/}h^{[S]}$. 
\end{definition}

We are now ready to define the  polynomial invariant.

\begin{definition}
Let D be an oriented diagram in the solid torus which represents a link L in the solid torus and let $w(D)$ be its writhe.

The Laurent polynomial $W_L(x,y,z,h) \in \mathbb{Z}[i][x,x^{-1},y,y^{-1},z,z^{-1},h]$ is defined as

$W_L = (xyz^{-1})^{-2w(D)}\sum_T \sum_S <D,T,S>$. 

Here the sums are over all Kauffman states T for the Alexander polynomial and all Kauffman states S for the Jones 
polynomial. 

It follows immediately from the definitions that $W_L(x,y,z,h)$ is homogenous in $x,y,z$ of degree 0. Therefore, we can 
replace it by $W_L = W_L(x,y,1,h)$ without loosing information.
\end{definition}

Let $D^2$ be the unit disc in the $(x,y)$-plane and let $V_1$ be the solid torus $D^2 \times \mathbb{R} \setminus z-axes$.
Let $V_2$ be the complementary solid torus $V \setminus V_1$.

\begin{definition}
A link L in the solid torus $V$ is called a {\em split link} if $L$ is isotopic to a link $ L_1 \cup L_2$ with the link $L_1$  
contained in $V_1$ and $L_2$ contained in $V_2$.  (If the axes $A$ is a meridian
of the link then our definition coincides with the usual definition of a split link in 3-space.)
\end{definition}

Let us consider oriented links in 3-space, i.e. links in the solid torus $V$ such that A is a positive meridian.
Then $W_L$ is linear in h and we can forget the variable h. 
Let us consider in this case the homogenous Laurent polynomial  $W_L(x,y,z)$ of
degree 0.

The following theorem is our first result.

\begin{theorem}
$W_L$ is an isotopy invariant for oriented links in the solid torus.

$W_L = 0$ for each split link L in the solid torus.

\end{theorem}

\begin{example}
Let K be the right-handed trefoil in 3-space. Then

$W_K = (x^4+x^2y^2+y^4)((xy)^{-10}+(xy)^{-8}+(xy)^{-4})$

The polynomial for the left-handed trefoil K! is obtained by replacing x and y by their inverses (which is not an immediate
corollary of the definitions).
\end{example}

I am very grateful to Stepan Orevkov, who has written a computer program in order to calculate $W_L$ (see \cite{O}). Calculations with this program 
suggest the following conjecture.

\begin{conjecture}
Let $L$ be a link in 3-space. Then $W_L$ is the product of a homogenous polynomial in $x$ and $y$ with a polynomial
in $xy$. The homogenous polynomial coincides for $y = ix^{-1}$ with the Alexander polynomial $\Delta_L(x^4)$.
\end{conjecture}

\begin{remark}
In order to find the  polynomial $W_L$, we have of course associated to each of the eight pictures in Fig. 3 a new variable 
as well as to the corresponding pictures for a negative crossing together with the variable $d$. Invariance under the 
Reidemeister moves has led to a 
non-linear system of 17 variables and 60 equations. It turns out that this system has a unique solution which gives a 
homogenous Laurent polynomial of degree 0 of three variables . 

I am very grateful to Benjamin Audoux and to Delphine Boucher for their help in solving the above system.
\end{remark}

\begin{remark}
We have developed a machinery, called {\em one parameter knot theory}, which produces new knot polynomials 
from state sums in the solid torus of classical link polynomials (see \cite{F3} and also \cite{FK}, which contains some 
necessary preparations). This has worked perfectly for the Kauffman state sums of the Alexander and of the Jones
polynomials (see \cite{F3}). 

The present paper is a result of our search for new state sums.The new state sum takes into account
not only the crossings and the arcs in the knot diagram which connect them, but also the components of its complement
 in the annulus. 

\end{remark}

State sums turned out to be very useful in order to categorify link polynomials in a combinatorial way 
(compare \cite{Kh} and \cite{MOS}).

\begin{question}
Can $W_L$ be categorified as a whole ?

\end{question}

We refine now $W_L$ with a new variable for a special class of links.

Let $L = L_1 \cup L_2 \cup L_3$ be an oriented link in 3-space such that $L_1 \cup L_2$ is a Hopf link and
$L_3$ is an arbitrary link.
We consider $L$ up to isotopy which preserves this decomposition. Consequently, instead of $L$ in $S^3$ we can 
equivalently study $L_3$ in the thickened torus $S^3 \setminus (L_1 \cup L_2) = T^2 \times \mathbb{R}$ (compare e.g.
\cite{F1}).  Moreover, the 
natural projection 
of $L_1$ and $L_2$ into the 2-torus $T^2$ determines a distinguished pair of generators , say $a$ and $b$, of 
$H_1(T^2)$.

Instead of the projection of $L$ into the annulus we use now the projection of $L_3$ into the 2-torus in order to define
$W(L_3 \hookrightarrow T^2 \times \mathbb{R})$ in exactly the same way as before besides the following two 
changings:

-there are no star regions, because the Euler characteristic of the torus vanishes (if there are no Kauffman states
for the Alexander polynomial, then the invariant vanishes)

-the non contractible circles of a double state are all parallel in $T^2$ and hence represent the same homology class
$\vert ma+nb\vert$. Here $m$ and $n$ are coprime integers and the homology class $ma+nb$ is only well defined
up to sign. Hence, we have to replace the variable $h$ by the homology class $\vert ma+nb\vert$. Contractible circles
are traded to factors $d$ as previously.

However, it turns out that $W(L_3 \hookrightarrow T^2 \times \mathbb{R})$ can be refined with a new variable, which 
comes from an unexpected relation between the Kauffman states for the Alexander polynomial and those for the
Jones polynomial !

\begin{definition}
A dot in a double state $(T,S)$ is called {\em counting} if the spitting at the corresponding crossing is not in the direction
of the dot. Each counting dot is {\em nearest} to a unique circle in the double state.
\end{definition}

For each circle in a double state we consider all its nearest counting dots on its left side and on its right side.

We show an example of a circle with exactly one nearest counting dot in Fig. 4.

\begin{figure}
\centering 
\psfig{file=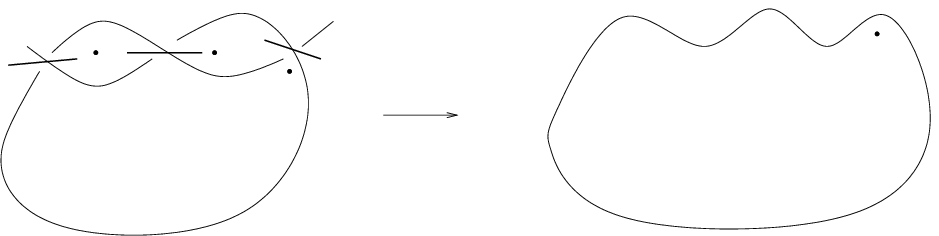}
\caption{}
\end{figure}

\begin{definition}
Let $C$ be a circle in a double state $(T,S)$. We chose its left side and its right side arbitrarily. Let $v_+$ be the number
of nearest counting dots on its right side and let $v_-$ be the number of nearest counting dots on its left side. The 
{\em weight}  of the circle is then the natural number $v(C) = \vert v_+ - v_-\vert$.
\end{definition}

\begin{remark}
We could define the weight $v(C)$ in exactly the same way in our previous case of links in the solid torus. 
Unfortunately, it turns out that each contractible circle has always the weight 1 and each non contractible circle has always
the weight 0.

But in the case of the thickened torus the weight can be arbitrarily large for non contractible circles.
We show an example in Fig. 5. The given double state has exactly  two non contractible circles and each has the 
weight two.
\end{remark}

\begin{remark}
If we choose an orientation on $C$ then the weight $v$ could be defined as an integer, because
the right side and left side are now defined canonically.
\end{remark}

\begin{figure}
\centering 
\psfig{file=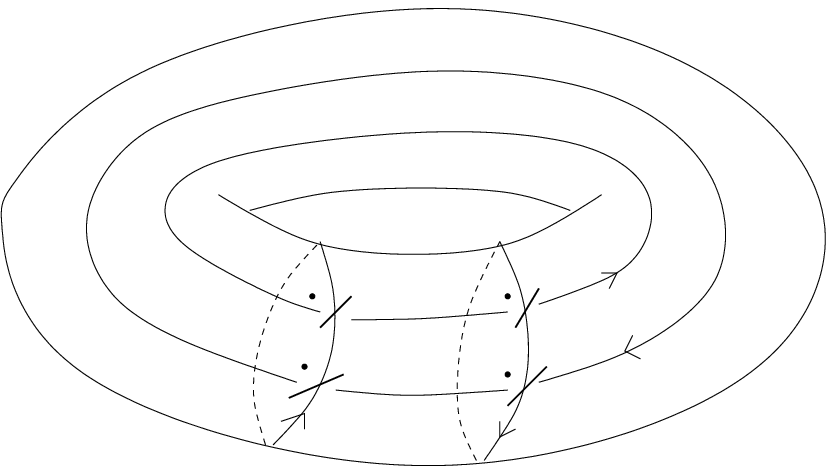}
\caption{}
\end{figure}

We make the following observation (and we left the verification to the reader) : all non contractible circles in a double state
have always the same weight and each contractible circle has the weight 1.
Consequently, instead of a weighted configuration of non contractible circles all the information is already given by 
the following data: let $C$ be a non contractible circle in the double state $(T,S)$. Then we define

- the number of (parallel) non contractible circles , denoted by $s(C)$

- the (absolut value) of its homology class (with respect to the distinguished set of generators), denoted by $\vert C\vert$

- the weight of a non contractible circle, encoded by  $t^{v(C)}$, where $t$ is a new variable.

\begin{definition}
Let $D$ be the diagram in $T^2 \times \mathbb{R}$ represented by $L_3$. Then the {\em refined} double bracket $<D,T,S>_H$ is the product of 
the previously
defined monomials over all crossings of $D$ and $d^{/S/}\vert C\vert^{s(C)}t^{v(C)}$. (Here, $\vert C\vert^{s(C)}$ is just the notation
for $s(C)$ parallel curves in the given homology class (up to sign) $\vert C\vert$.)
\end{definition}

\begin{definition}
The {\em refined} polynomial $W^H_L(x,y,z,t,\vert C\vert)$ is defined as

$W^H_L = (xyz^{-1})^{-2w(D)}\sum_T \sum_S <D,T,S>_H$. 

Here $<D,T,S>_H$ is the refined double bracket.

\end{definition}

$W^H_L$  can be seen as a 3-variable polynomial for each fixed element in $H_1(T^2)$.
 
\begin{theorem}
$W^H_L$ is an isotopy invariant for those oriented links in the 3-sphere which contain a distinguished  Hopf link as a sublink.
\end{theorem}

\begin{example}
Let $L$ be the link shown in Fig. 6. There is only one crossing and hence there are  only eight double states.
One easily calculates 

$W^H_L = (x^2+y^2-2ixyt)\vert a+b\vert + (2+ixy^{-1}t+ix^{-1}yt)\vert a-b\vert$

(we have set $z=1$).
\end{example}

\begin{figure}
\centering 
\psfig{file=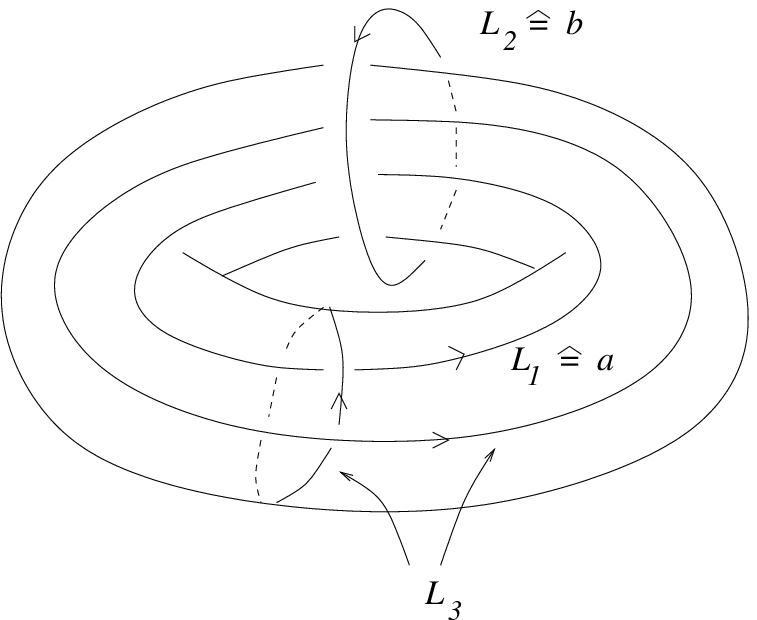}
\caption{}
\end{figure}

\begin{question}
Can $W^H_L$ detect non invertibility of $L$ ?

Can it detect mutations which preserve the Hopf link $L_1 \cup L_2$ ?
\end{question}

\begin{remark}
$W^H_L$ could be easily generalized for links in handle bodies, but it would no longer lead to an invariant for classical
 links.
\end{remark}

\section{Proofs}

We have to check that $W_L$ is invariant under the Reidemeister moves (compare e.g. \cite{BZ}).

{\em Reidemeister I \/}

The four Reidemeister I moves are shown in Fig. 7.

\begin{figure}
\centering 
\psfig{file=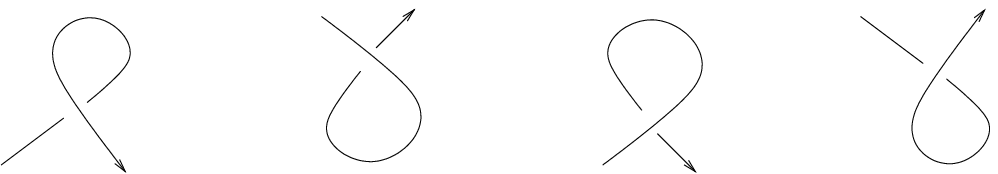}
\caption{}
\end{figure}

Notice, that the dot has to be in the newly created disc region, which can not be a star region. One easily calculates that
 the first two moves multiply
$\sum_T \sum_S <D,T,S>$ by $x^2y^2z^{-2}$  and the other two moves multiply it by $x^{-2}y^{-2}z^2$  and, consequently,
$W_L = (xyz^{-1})^{-2w(D)}\sum_T \sum_S <D,T,S>$ is invariant.

{\em Reidemeister II \/}

Let $IIa$ be the move where the two tangencies at the autotangency point have the same orientation and let $IIb$ be the
move where they have opposite orientations. Let us consider $IIa$ (the case $IIb$ as well as the mirror images of 
these moves are completely analogous and are left to
the reader). We splitt the set of all Kauffman states T for the Alexander
polynomial into four subsets shown in Fig. 8 up to Fig. 11. We assume here that there are no stars in the regions shown 
in the figures. In each of the figures we consider now the four Kauffman states S for the Jones polynomial.

\begin{figure}
\centering 
\psfig{file=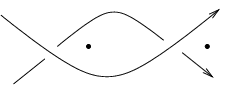}
\caption{}
\end{figure}

\begin{figure}
\centering 
\psfig{file=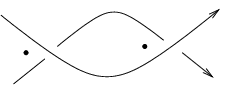}
\caption{}
\end{figure}

\begin{figure}
\centering 
\psfig{file=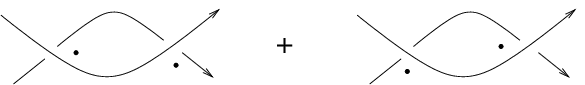}
\caption{}
\end{figure}

\begin{figure}
\centering 
\psfig{file=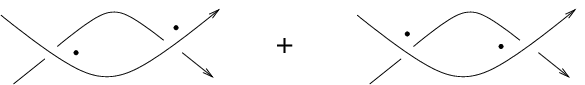}
\caption{}
\end{figure}

The results for two of the states are shown in Fig. 12 and in Fig. 13 (the remaining two states lead to identical pictures 
besides the dot, which has slided in its region to the right, respectively downwards). We omit to write brackets and a dot which is not in a quadrant of a crossing
means that there is a dot in the corresponding part of the region.

\begin{figure}
\centering 
\psfig{file=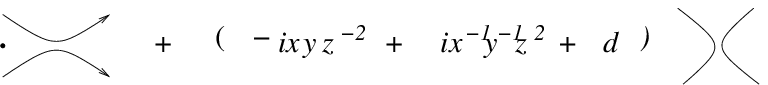}
\caption{}
\end{figure}

\begin{figure}
\centering 
\psfig{file=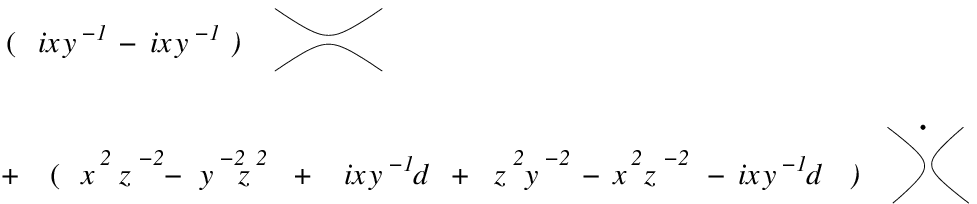}
\caption{}
\end{figure}

 It follows that $W_L$ is invariant under these Reidemeister II moves.
This is still true if there are stars in the regions, because again there can not be a star in the newly created region bounded
by the bigon. If the stars are situated as shown in Fig. 14 then it follows that $W_L = 0$. But notice that in this case $L$
 is necessarily a split link in the solid torus.

\begin{figure}
\centering 
\psfig{file=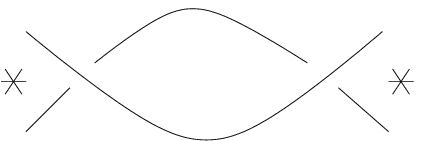}
\caption{}
\end{figure}

{\em Reidemeister III \/}

If one considers oriented diagrams 
then there are exactly eight different (local) types of Reidemeister III moves.
Let us call the {\em positive Reidemeister III move} those in which all three involved crossings are positive.
Fortunately, it turns out that in order to check that a polynomial is a knot invariant it suffices to check invariance only 
under  the positive Reidemeister III move and under all types of Reidemeister II moves (see e.g. Section 1 in \cite{F1} 
and also Sections 2.3 and 2.4 in \cite{F3}).

Under a Reidemeister III move a triangle component of the complement of the diagram shrinks to a point and the link 
diagram has an ordinary triple point in the projection. There are
three types of Kauffman states T for the Alexander polynomial at the triple point shown in Fig. 15.

\begin{figure}
\centering 
\psfig{file=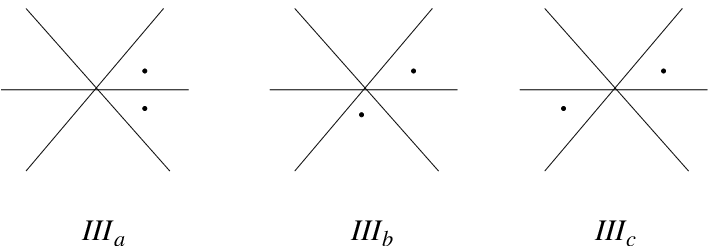}
\caption{}
\end{figure}

For each of the three 
types we will consider just
one of the six (respectively three) cases. The remaining cases are always analogue and are left to the reader. In the 
figures we draw only the planar image of the projection because all crossings are positive . (The rest of the states T
 are of course identical outside of the corresponding figures.)

{\em case IIIa}

We show the corresponding states T before and after the move in Fig. 16.

\begin{figure}
\centering 
\psfig{file=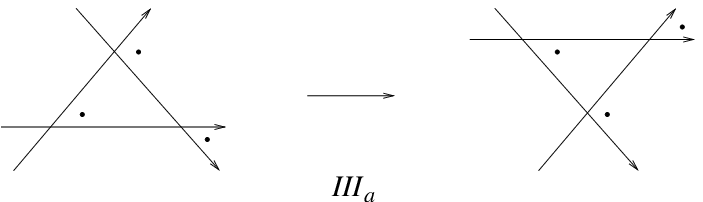}
\caption{}
\end{figure}

Fig. 17 shows now the contribution to $W_L$ before the move and Fig. 18 shows the contribution after the move.

\begin{figure}
\centering 
\psfig{file=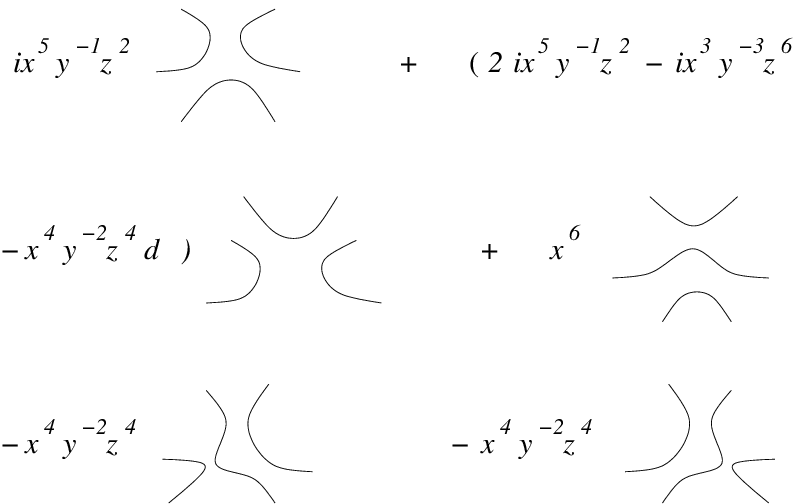}
\caption{}
\end{figure}

\begin{figure}
\centering 
\psfig{file=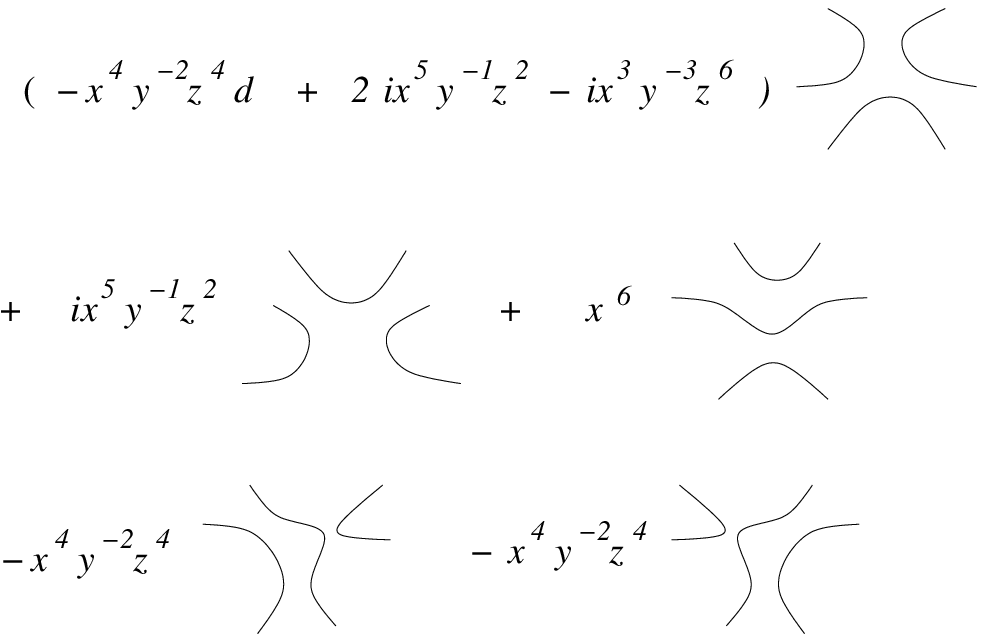}
\caption{}
\end{figure}

{\em case IIIb}

We show the corresponding states T before and after the move in Fig. 19.

\begin{figure}
\centering 
\psfig{file=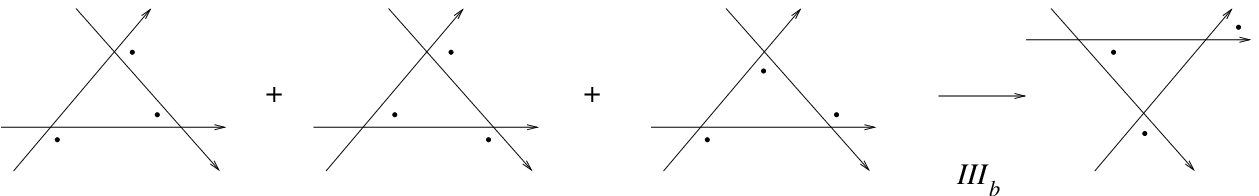}
\caption{}
\end{figure}

Fig. 20 shows the contribution to $W_L$ before the move and Fig. 21 shows the contribution after the move.

\begin{figure}
\centering 
\psfig{file=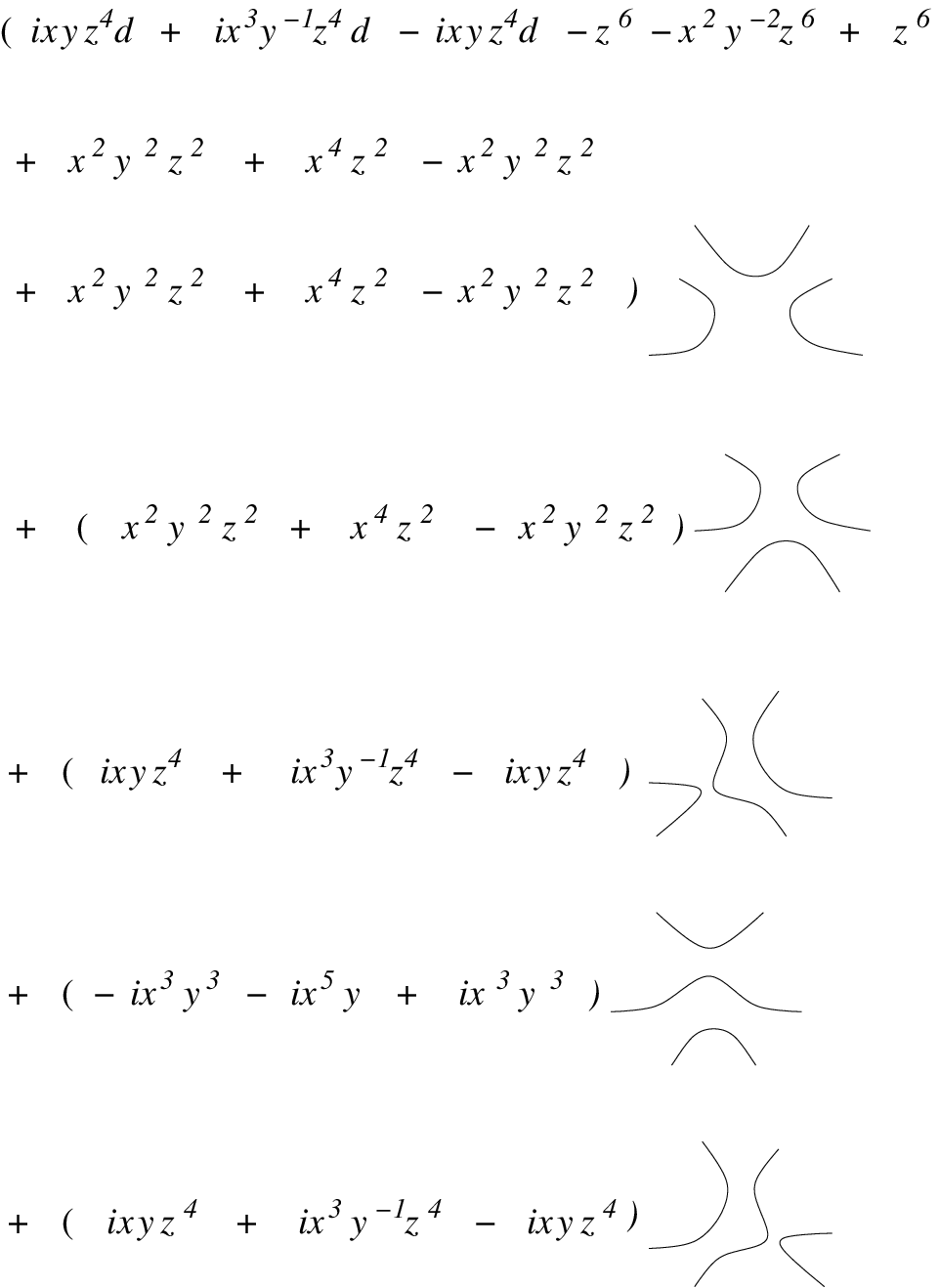}
\caption{}
\end{figure}

\begin{figure}
\centering 
\psfig{file=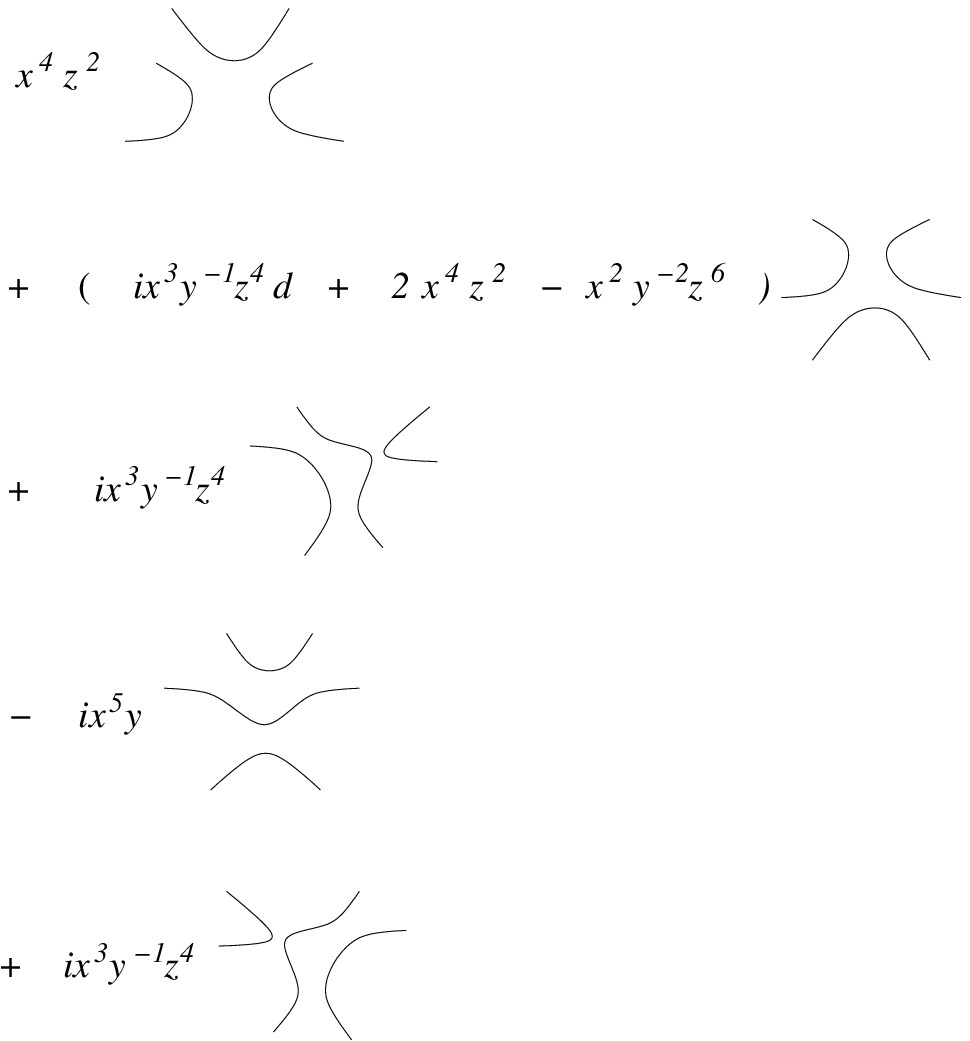}
\caption{}
\end{figure}

{\em case IIIc}

We show the corresponding states T before and after the move in Fig. 22.

\begin{figure}
\centering 
\psfig{file=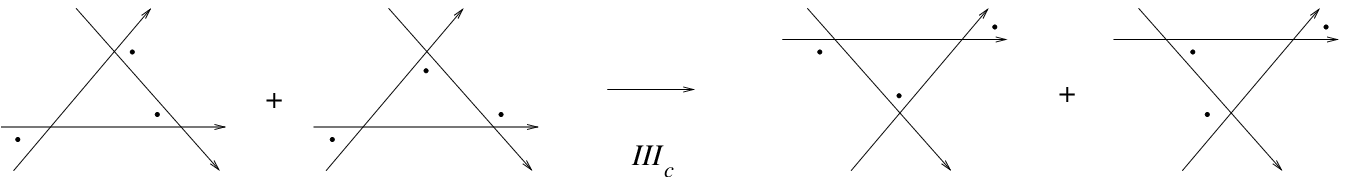}
\caption{}
\end{figure}

Fig. 23 shows the contribution to $W_L$ before the move and Fig. 24 shows the contribution after the move.

\begin{figure}
\centering 
\psfig{file=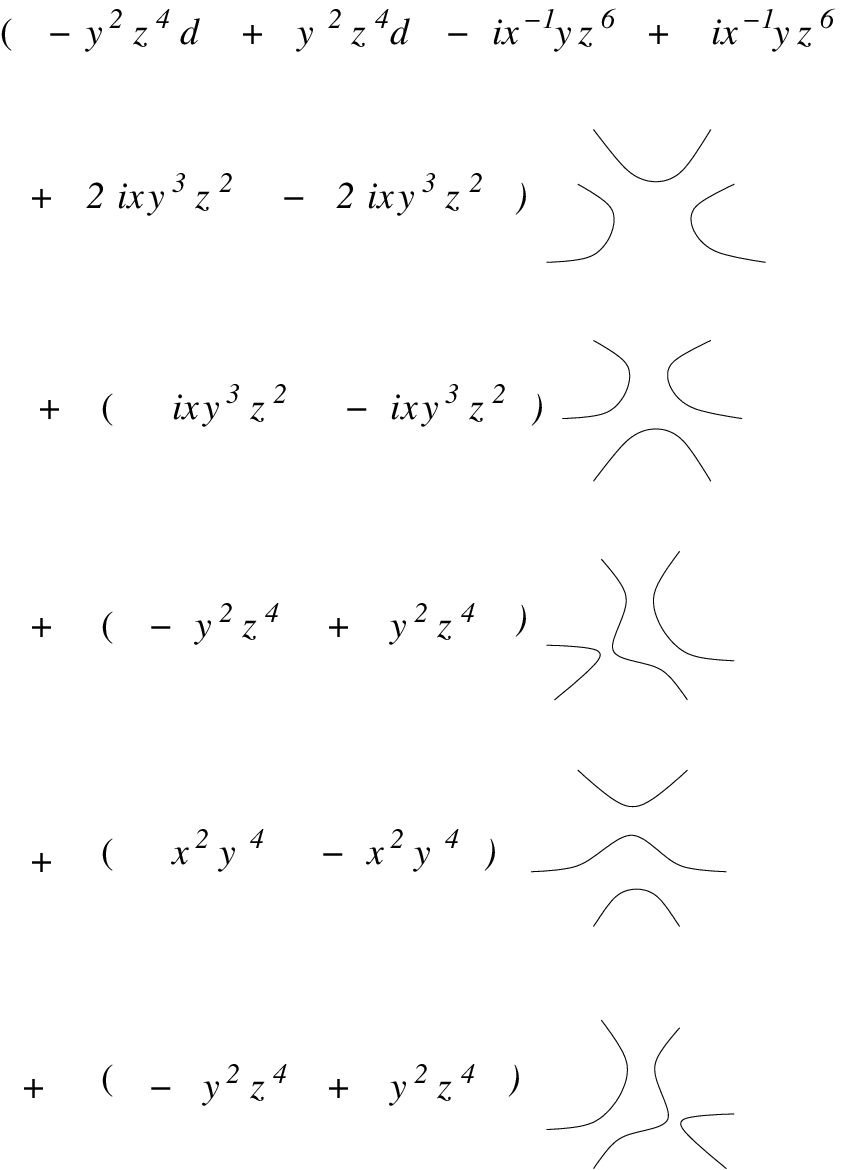}
\caption{}
\end{figure}

\begin{figure}
\centering 
\psfig{file=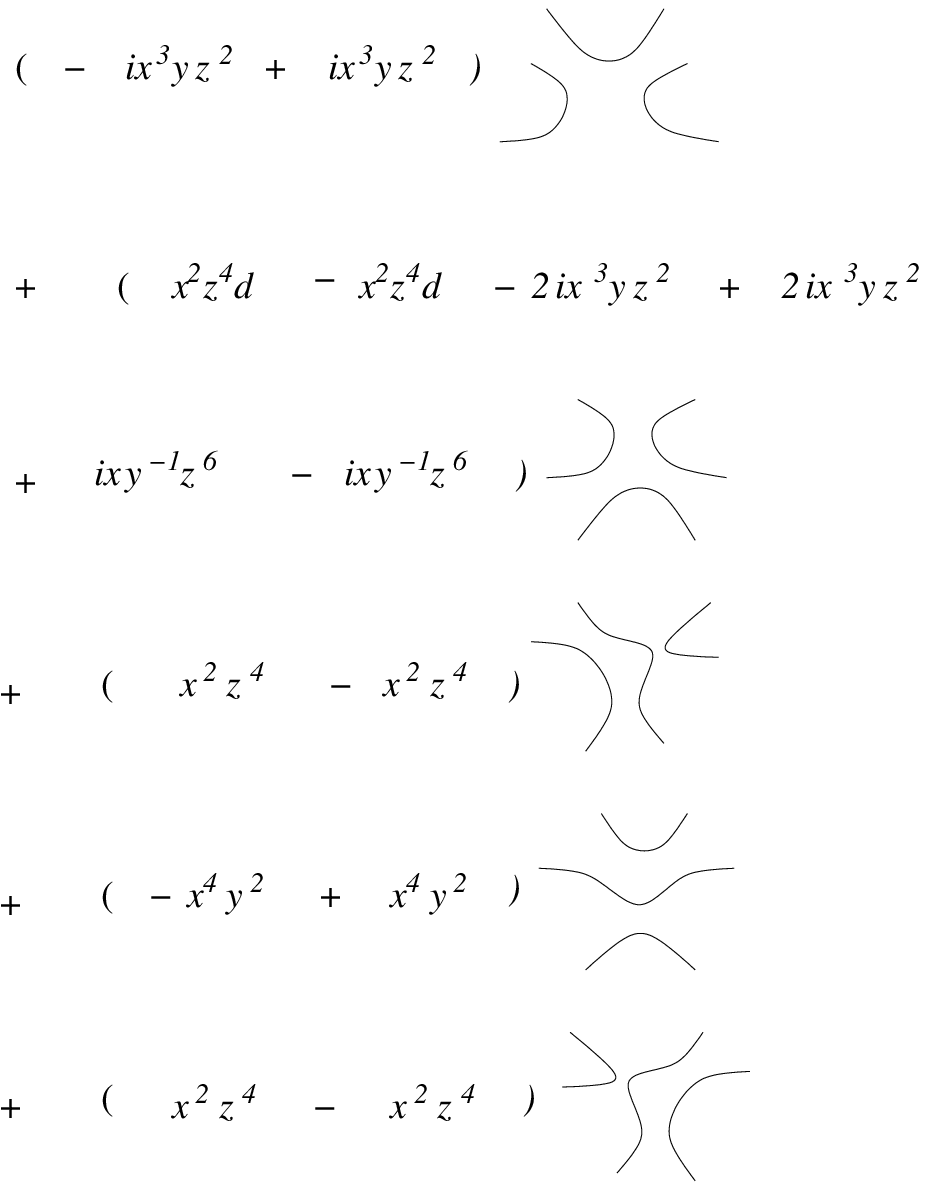}
\caption{}
\end{figure}
\begin{figure}
\centering 
\psfig{file=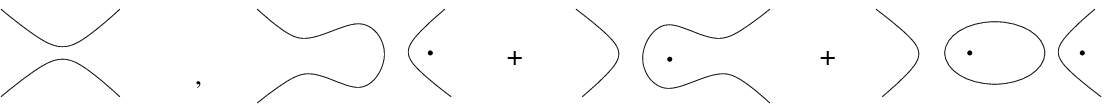}
\caption{}
\end{figure}

\begin{figure}
\centering 
\psfig{file=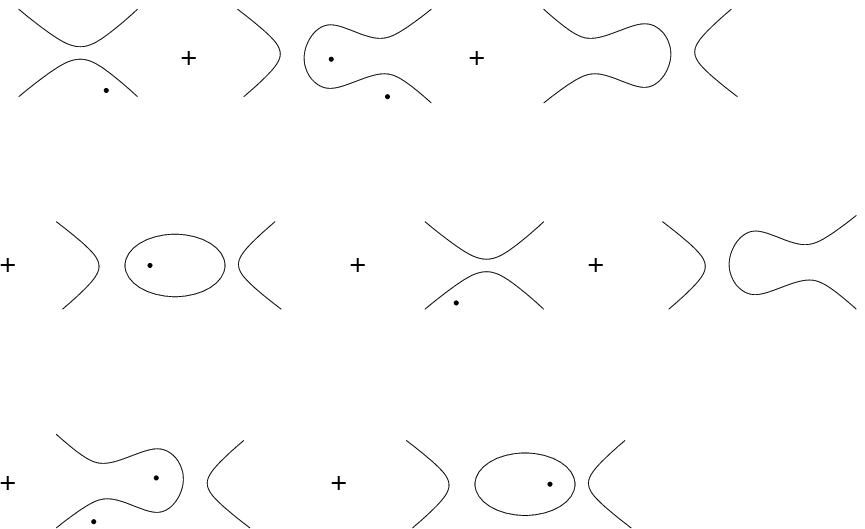}
\caption{}
\end{figure}
\begin{figure}
\centering 
\psfig{file=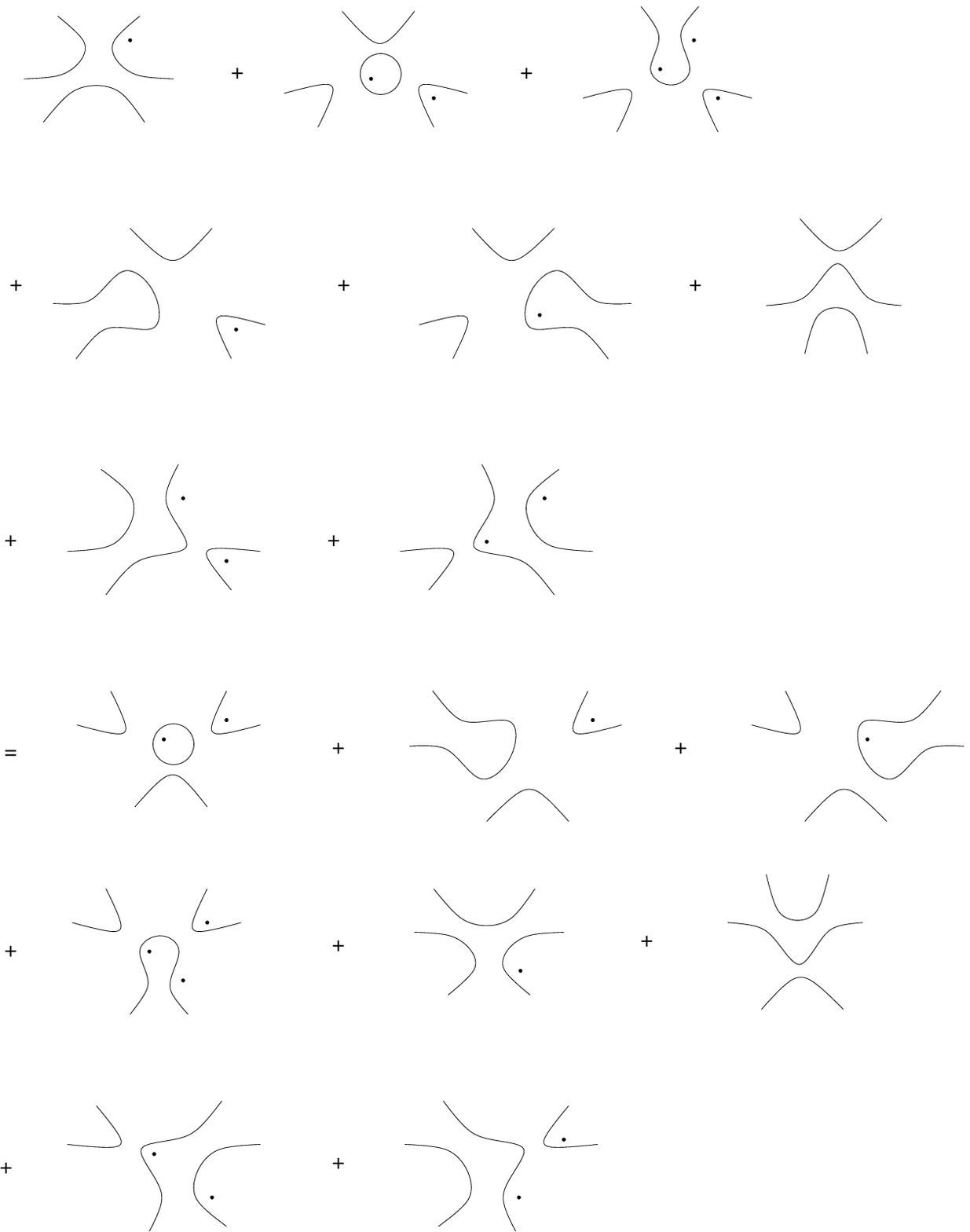}
\caption{}
\end{figure}

It follows that $W_L$ is invariant under these Reidemeister III moves.  Again, stars in the figures would not change the above
identities, because there can never be stars in the vanishing triangles.

Notice that the states of IIIc do not contribute at all to $W_L$. This was already observed for the Alexander polynomial 
in \cite{F3} , Remark 11.

Theorem 1 is proven.

\vspace{1cm}

The proof of Theorem 2 is completely analogous besides the following additional consideration:
we have to prove that all corresponding circles of configurations which enter into the same equation have the same 
weight $v$.
Notice that the weights are completely determined by the unoriented curves together with the dots in the torus $T^2$.
Therefore, it suffices to consider unoriented immersed curves in $T^2$ instead of oriented link diagrams.

{\em Reidemeister I \/}

The new dot is either not a counting dot or it is a counting dot nearest to a contractible circle.

{\em Reidemeister II \/}

Fig. 8 (as well as Fig. 9) leads to the configurations shown in Fig. 25. We draw only the counting dots.

Fig. 10 (as well as Fig. 11) leads to the configurations shown in Fig. 26. 
It follows that all corresponding circles have the same weight.

{\em Reidemeister III \/}

We consider as an example the case IIIa. The other cases are analogous and we left the verification to the reader.
Fig. 16 leads to the configurations shown in Fig. 27.
It follows again that all corresponding circles have the same weight. 

Theorem 2 is proven.

\vspace{1cm}

{\em Acknowledgements}--- I wish to tank Benjamin Audoux, Delphine Boucher and Stepan Orevkov for their help.

Laboratoire de Math\'ematiques

Emile Picard

Universit\'e Paul Sabatier

118 ,route de Narbonne 

31062 Toulouse Cedex 09, France

fiedler@picard.ups-tlse.fr

\end{document}